\documentclass{elsart}

\usepackage{amsfonts}
\usepackage{amsmath}
\usepackage{amssymb}
\input{xy}
\xyoption{all}

\newcommand{\DA}{\mathcal{D}}

\newcommand{\rad}{\textup{rad}\,}
\newcommand{\Hom}{\textup{Hom}\,}

\newcommand{\Ext}{\textup{Ext}}
\newcommand{\Ann}{\textup{Ann}}
\newcommand{\End}{\textup{End}}

\newcommand{\ind}{\textup{ind}\,}
\newcommand{\add}{\textup{add}\,}

\newcommand{\pd}{\textup{pd}\,}
\newcommand{\isomorphe}{\cong}
\newcommand{\CA}{{\mathcal{C}_A}}

\newcommand{\ot}{\leftarrow}

\newcommand{\zd}{\delta}
\newcommand{\zG}{\Gamma}
\newcommand{\zg}{\gamma}
\newcommand{\zS}{\Sigma}

\newcommand{\za}{\alpha}
\newcommand{\zb}{\beta}

\newcommand{\ze}{\epsilon}

\addtocounter{section}{-1}
\begin{document}

\begin{frontmatter}

\title{Cluster-tilted algebras and slices}
  \author{Ibrahim Assem\thanksref{nserc}}, \address{D\'epartement de Math\'ematiques,
  Universit\'e de Sherbrooke, Sherbrooke (Qu\'ebec), J1K 2R1, Canada}
  \ead{ibrahim.assem@usherbrooke.ca}
  \author{Thomas Br\"ustle\thanksref{nserc+}}, \address{D\'epartement de Math\'ematiques,
  Universit\'e de Sherbrooke, Sherbrooke (Qu\'ebec), J1K 2R1, Canada \emph{and} 
  Department of Mathematics, Bishop's University,  Lennoxville, (Qu\'ebec),
  J1M 1Z7, Canada} 
  \ead{thomas.brustle@usherbrooke.ca}
  \author{Ralf Schiffler\thanksref{umass}}, \address{Department of Mathematics and
  Statistics, University of Massachusetts at Amherst, Amherst, MA
  01003-9305, USA}
  \ead{schiffler@math.umass.edu} 

  \thanks[nserc]{Partially supported by the NSERC of Canada and the University of Sherbrooke}
  \thanks[nserc+]{Partially supported by the NSERC of Canada and the universities of Sherbrooke and Bishop's}
  \thanks[umass]{Partially supported by the University of Massachusetts}

\begin{abstract}
We give a criterion allowing to verify whether or not two tilted
algebras have the same relation-extension (thus correspond to the same
cluster-tilted algebra). This criterion is in terms of a combinatorial
configuration in the Auslander-Reiten quiver of the cluster-tilted
algebra, which we call local slice.
\end{abstract}
\begin{keyword}
Local slice \sep cluster-tilted algebra \sep relation-extensions of  tilted algebras
\MSC 16S70 \sep 16G20 
\end{keyword}
\end{frontmatter}


\begin{section}{Introduction}\label{sect intro}
Cluster categories were introduced in \cite{BMRRT} and, for type $\mathbb{A}$,
also in \cite{CCS1}, as a categorical model allowing to understand
better the cluster algebras of Fomin and Zelevinsky
\cite{FZ1}. Cluster-tilted algebras were defined in \cite{CCS1}
for type $\mathbb{A}$, and in \cite{BMR} for arbitrary hereditary algebras as
follows:
Let $A$ be a hereditary algebra and $\tilde T$ be a tilting object in
the associated cluster category $\CA$, that is, an object such that
$\Ext^1_\CA(\tilde T,\tilde T)=0$ and the number of isomorphism
classes of indecomposable summands of $\tilde T$ equals the rank of
the Grothendieck group of $A$, then the algebra $B=\End_\CA\tilde T$ is called cluster-tilted. These algebras have been
studied by several authors, see, for instance,
\cite{ABS,BMR,BMR3,BR,CCS2,KR}. In particular, they were shown in
\cite{ABS} to be closely related to the tilted algebras introduced by
Happel and Ringel  in the early eighties \cite{HR}. Indeed, let $C$ be
a tilted algebra, then the trivial extension $\tilde C= C\ltimes
\Ext^2_C(DC,C)$ of $C$ by the $C$-$C$-bimodule $\Ext^2_C(DC,C) $ is
cluster-tilted, and every cluster-tilted algebra is of this
form. Thus, we have a surjective map $C\mapsto \tilde C$ from tilted
to cluster-tilted algebras. However, easy examples show that this map
is not injective. Our objective in this paper is to give a criterion
allowing to verify whether for two tilted algebras $C_1$ and
$C_2$, the corresponding cluster-tilted algebras $\tilde C_1$ and
$\tilde C_2$ are isomorphic or not.

Since tilted algebras are characterised by the existence of complete
slices in their Auslander-Reiten quiver (see, for instance,
\cite{HR,R,L2,Sk1} or \cite{ASS}), it is natural to study the
corresponding concept for cluster-tilted algebras. For this purpose,
we introduce what we call a \emph{local slice}, by weakening the
axioms of complete slice (thus, in a tilted algebra, complete slices
are local slices). We show that a complete slice in a tilted algebra
$C$ embeds as a local slice in $\tilde C$ (and, in fact, any local
slice in $\tilde C$ is of this form).

Our main theorem is the following:

\begin{thm}
Let $B$ be a cluster-tilted algebra. Then a tilted algebra
$C$ is such that $B=C\ltimes \Ext^2_C(DC,C) $ if and only if
there exists a local slice $\zS$ in $\textup{mod}\, B$ such
that \[C=B/\Ann_{B}\zS.\]
\end{thm}

We also show that cluster-tilted algebras have many local slices. In
fact, all but at most finitely many indecomposable modules lying in
the transjective component of the Auslander-Reiten quiver belong to a
local slice. If the cluster-tilted algebra is of tree type (which is
the case, for instance, if it is of Dynkin type or of Euclidean type
distinct from $\tilde{ \mathbb{A}}$), then this is the case for all
indecomposable modules in this component.

We now describe the contents of our paper.
In the first section, we introduce the notion of local section in a
translation quiver and in the second section we study sections and local
sections in the derived category. In the third section, 
we introduce the concept of a local slice and prove our main result,
and in section four, we prove that cluster-tilted algebras of tree
type have sufficiently many local slices.

We would like to thank David Smith for useful comments on the paper.
\end{section}

\begin{section}{Preliminaries on translation quivers}\label{sect 1}

\begin{subsection}{Notation} Throughout this paper, all algebras are
  connected finite dimensional algebras over an algebraically closed
  field $k$. For an algebra $C$, we denote by $\textup{mod}\, C$ the
  category of finitely generated right $C$-modules and by $\ind C$ a
  full subcategory of $\textup{mod}\, C$ consisting of exactly one
  representative from each isomorphism class of indecomposable
  modules. When we speak about a $C$-module (or an indecomposable
  $C$-module), we always mean implicitly that it belongs to
  $\textup{mod}\, C$ (or to $\ind C$, respectively). Also, all
  subcategories of $\textup{mod}\, C$ are full and so are identified
  with their object classes. Given a subcategory $\mathcal{C}$ of
  $\textup{mod}\, C$, we sometimes write $M\in \mathcal{C}$ to express
  that $M$ is an object in $\mathcal{C}$. We denote by $\add
  \mathcal{C}$ the full subcategory of $\textup{mod}\, C$ with objects
  the finite direct sums of modules in $\mathcal{C}$ and, if $M$ is a
  module, we abbreviate $\add\{M\}$ as $\add M$. We denote the
  projective (or injective) dimension of a module $M$ as $\pd M$ (or
  $\textup{id}\, M$, respectively). The global dimension of $C$ is denoted by
  gl.dim.$C$ and its Grothendieck group by $K_0(C)$. Finally, we
  denote by $\zG(\textup{mod}\, C)$ the Auslander-Reiten quiver of an
  algebra $C$, and by $\tau_C= D\, Tr$, $\tau^{-1}_C=Tr\, D$ its
  Auslander-Reiten translations. For further definitions and facts
  needed on $\textup{mod}\, C$ or $\zG(\textup{mod}\, C)$, we refer
  the reader to  \cite{ASS}. We also need facts on the bounded
  derived category $\mathcal{D}^b(\textup{mod}\, C)$ of
  $\textup{mod}\, C$, for which we refer to \cite{H}.
\end{subsection} 

\begin{subsection}{Sections}\label{sect 2.1}

For translation quivers, we refer to \cite{ASS,R}. Let 
 $(\zG,\tau)$ be a connected translation
quiver. We recall that a path $x=x_0\to x_1\to \ldots \to x_t=y$ in
$\zG$ is called \emph{sectional} if, for each $i$ with $0<i<t$, we
have $\tau\,x_{i+1}\ne x_{i-1}$. A full connected subquiver $\zS$ of
$\zG$ is said to be \emph{convex} in $\zG$ if, for any path $x=x_0\to
x_1\to \ldots \to x_t=y$ in $\zG$ with $x,y\in \zS_0$, we have $x_i\in
\zS_0$ for all $i$. It is called \emph{acyclic}  if there is no cycle
$x=x_0\to x_1\to \ldots \to x_t=x$ (with $t>0$) which is entirely contained in
$\zS$. The following definition is due to Liu and Skowro\'nski (see
 \cite{L1,Sk1} or else \cite{ASS}).

\begin{defn}
Let $(\zG,\tau)$ be a connected translation quiver.
A connected full subquiver  $\zS$  of $\zG$ is a \emph{section} in
$\zG$ if:
\begin{itemize}
\item[(S1)] $\zS$ is acyclic.
\item[(S2)] For each $x\in \zG_0$, there exists a unique $n\in
  \mathbf{Z}$ such that $\tau^n\,x\in \zS_0$.
\item[(S3)]  $\zS$ is convex in $\zG$.
\end{itemize}   
\end{defn}     

This definition is  motivated by the study of tilted algebras. The
  well-known criterion of Liu and Skowro\'nski asserts that, if $C$ is an
  algebra, and $\zS$ is a faithful section in a component of its
  Auslander-Reiten quiver such that $\Hom_C(X,\tau_CY)=0$ for all
  $X,Y\in \zS_0$, then $C$ is tilted having $\zS$ as complete slice
  (see \cite{L1,Sk1,ASS}).

We note that, if a translation quiver $\zG$ contains a section, then
$\zG$ is acyclic. 

\end{subsection}

\begin{subsection}{Presections}\label{sect 2.2}
We need some weaker notions. The first one is the following.

\begin{defn}
Let $(\zG,\tau)$ be a connected translation quiver.
A connected full subquiver  $\zS$  of $\zG$ is called a \emph{presection} in
$\zG$ if it satisfies the following two conditions:
\begin{itemize}
\item[(P1)] If $x\in \zS_0$ and $x\to y$ is an arrow, then either
  $y\in \zS_0$ or $\tau y\in \zS_0$.
\item[(P2)] If $y\in \zS_0$ and $x\to y$ is an arrow, then either
  $x\in \zS_0$ or $\tau^{-1}x\in \zS_0$.
\end{itemize}   
\end{defn}     

The following Lemma collects the elementary properties of presections.
 Recall that a translation quiver is called \emph{stable} if there are
 neither projective,
 nor injective points in $\zG$.

\begin{lem}\label{lem 2.2} Let $(\zG,\tau)$ be a connected translation quiver.
\begin{itemize}
\item[(a)] If $\zS$ is a  section in $\zG$, then $\zS$  is a presection.
\item[(b)] Any path entirely contained in a presection $\zS$ is a
  sectional path. 
\item[(c)] If the translation quiver $\zG$ is stable,
   then conditions (P1) and (P2) are equivalent.
\item[(d)]  If $\zS$ is a  presection in a stable translation quiver
  $\zG$, then $\zS$   intersects every $\tau$-orbit of $\zG$ at least once.
\end{itemize}   
\end{lem}  

\begin{pf}
(a) This is well-known (see, for instance, \cite[VIII.1.4 p.304]{ASS}).

(b) Assume that $\zS$ is a presection in $\zG$, and that $x=x_0\to
  x_1\to \ldots \to x_t=y$ is a path lying entirely in $\zS$. If this
  path is not sectional, then there exists a least $i$ with $0<i<t$
  and $\tau\,x_{i+1}=x_{i-1}.$ But, in this case, we have arrows
  $x_{i-1}\to x_i$ and $x_i\to \tau^{-1}\,x_{i-1}=x_{i+1}$ with
  $x_i,x_{i-1},x_{i+1} \in \zS_0$, a contradiction.

(c) Assume (P1) holds and that $x\to y$ is an arrow with $y \in
  \zS_0$. Since $x$ is not injective, there exists an arrow $y\to
  \tau^{-1}x$. Applying (P1) to the latter yields that $\tau^{-1}x\in
  \zS_0$ or $x=\tau(\tau^{-1}x)\in \zS_0$. Thus (P2) holds. The
  converse is shown in the same way.

(d) It suffices to prove that, if $x\in \zS_0$ and $y\in\zG_0$ lie in
  two neighbouring $\tau$-orbits, then $\zS$ intersects the
  $\tau$-orbit of $y$. Since $\zG$ is stable, there exists $m\in
  \mathbf{Z}$ such that there is an arrow $x\to \tau^m y$. But then
  $\tau^m y\in \zS_0$ or $\tau^{m+1}y\in \zS_0$.
\qed
\end{pf}  
 
For instance, let $\zG$ be a stable tube. A ray in $\zG$ is a  presection
  but clearly not a section. On the other hand, if $\zG$ is  a
  component of type $\mathbf{Z}A_\infty$, then a ray is a section.

\end{subsection}

\begin{subsection}{Local sections}\label{sect 2.3}
We now define an intermediate notion between those of presection and section.

\begin{defn}
Let $(\zG,\tau)$ be a connected translation quiver. A presection $\zS$
in $\zG$ is called a \emph{local section} if it satisfies moreover the
following additional condition:
\begin{itemize}
\item[] $\zS$  is sectionally convex, that is, if $x=x_0\to x_1\to
  \ldots \to x_t=  y$ is a sectional path in $\zG$ such that $x,y\in
  \zS_0$, then  $x_i\in\zS_0$ for all $i$.
\end{itemize}   
\end{defn}     

The following Lemma is immediate.

\begin{lem}\label{lem 2.3}
Any section is a local section.
\end{lem}  

\begin{pf}
Indeed, any section is convex, and hence sectionally convex. We then
apply Lemma \ref{lem 2.2} (b).
\qed
\end{pf}  
\end{subsection}

The main result of this section is the following.

\begin{prop}\label{prop 2.4}
Let $Q$ be a finite acyclic quiver and $\zG=\mathbf{Z}Q$. The
following conditions are equivalent for a connected full subquiver
$\zS$ of $\zG$ such that $|\zS_0|=|Q_0|$:
\begin{itemize}
\item[(a)] $\zS $ is a section.
\item[(b)] $\zS $ is a local section.
\item[(c)] $\zS$ is a presection.
\end{itemize}   
\end{prop}

\begin{pf}
Because of Lemma \ref{lem 2.3} (and the definition of local section),
it suffices to prove that (c) implies (a). Since $\zG$ is a stable
translation quiver, then, because of Lemma \ref{lem 2.2} (d), $\zS$
intersects every $\tau$-orbit of $\zG$ at least once. Then, it follows
from the hypothesis that $|\zS_0|=|Q_0|$ that $\zS$ intersects each
$\tau$-orbit of $\zG$ exactly once. Since $\zS$ is clearly acyclic
(because $\zG$ is), there only remains to prove convexity. Suppose
that there exists a path $x=x_0\to x_1 \to \ldots \to x_t=y$ in $\zG$
with $x,y \in \zS_0$ but $x_1,\ldots, x_{t-1}\notin \zS_0$ $(t\ge
2)$. Since $\zS$ is a presection, we have  $\tau\,x_1\in \zS_0$. From
the arrow $\tau\,x_1\to \tau\,x_2$, we deduce that either
$\tau\,x_2\in \zS_0$ or $\tau^2x_2\in\zS_0$. Repeating this argument
$t$ times, we get that one of $\tau\, y, \tau^2 y, \ldots, \tau^t y$
lies in $\zS$. This contradicts the fact that $\zS$ cuts each
$\tau$-orbit exactly once, because $y\in \zS_0$.
\qed
\end{pf}

\end{section}


\begin{section}{Sections and local sections in the derived
    category.}\label{sect 3} 

Let $A$ denote a hereditary algebra, and $\DA=\mathcal{D}^b(\textup{mod}\,A)$ denote the
bounded derived category of $\textup{mod}\, A$. We denote by
$\tau_{\DA}, \tau_{\DA}^{-1}$ the Auslander-Reiten
translations in $\DA$. We recall that the Auslander-Reiten
quiver of $\DA$ consists of two types of components: the
regular (which correspond to the regular components of
$\zG(\textup{mod}\, A)$ and their shifts) and the transjective (which
are the form $\mathbf{Z}Q$, where $Q$ denotes the ordinary quiver of
$A$), see \cite{H}. 

We need the following result, known as
Skowro\'nski's Lemma \cite{Sk2,ASS}.

\begin{lem}\label{lem 3.1}
Let $C$ be an artin algebra. Assume that a $C$-module $M$ is the
direct sum of $m$ pairwise non-isomorphic indecomposable modules and
is such that $\Hom_C(M,\tau_C M)=0$. Then $m\le \textup{rk K}_0(C)$.
\end{lem}

The following Lemma is motivated by \cite{Sk1,L2}.

\begin{lem}\label{lem 3.2}
Let $\zS$ be a section in a connected component $\zG$ of
$\zG(\mathcal{D})$ then the following conditions are equivalent: 
\begin{itemize}
\item[(a)] $\zS$ is convex in $\ind \DA$, that is, 
if $X=X_0\to X_1\to \ldots \to
  Y$ is a sequence of non-zero morphisms between indecomposable
  objects in  $\DA $ such that  $X,Y\in\zS_0$, then
  $X_i\in\zS_0$ for all $i$.
\item[(b)] $\Hom(X,\tau_{\DA} Y)=0 $ for all $X,Y \in\zS_0$.
\item[(c)] $\Hom(\tau^{-1}_{\DA} X,Y)=0 $ for all $X,Y \in\zS_0$.
\item[(d)] $\zS_0$ is finite.
\item[(e)] $|\zS_0|=\textup{rk K}_0(A)$.
\item[(f)] $\zG$ is a  transjective component.
\end{itemize}    
\end{lem}  

\begin{pf}
The equivalence of  $(b)$ and $(c)$ follows trivially from the fact that $\tau$ is an automorphism of $\DA$.

$(a)$ implies $(b)$. A non-zero morphism $X\to \tau_{\DA}\, Y$ induces
a path $X\to \tau_{\DA}\, Y \to * \to Y$ in $\ind \DA$. The hypothesis
implies that both $Y$ and $\tau_{\DA}\,Y$ lie in $\zS$, a contradiction. 

$(b)$ implies $(f)$. If $\zG$ is not a transjective component, then
we may assume without loss of generality that $\zG$ is concentrated in
degree zero. By Skowro\'nski's Lemma \ref{lem 3.1}, $\zS$ is finite. Now
this contradicts the fact that $\zS$ intersects each
$\tau_{\DA}$-orbit exactly once.

$(f)$ implies $(e)$. The number of $\tau_{\DA}$-orbits in a
transjective component is exactly $\textup{rk K}_0(A)$.

$(e)$ implies $(d)$. This is trivial.

$(d)$ implies $(a)$. Suppose that $X=X_0\stackrel{f_1}{\to}
X_1\stackrel{f_2}{\to}  \ldots\stackrel{f_t}{\to} X_t=Y$ is a path,
where the $f_i$ are non-zero morphisms, the $X_i$ lie in $\zG$ and 
$X,Y\in \zS_0$. 
Now, since $\zS_0$ is finite, then $\zG$ has finitely many
$\tau_\DA$-orbits. Therefore, $\zG $ is transjective. This implies that 
the $f_i$ lie in a finite power of the radical of
$\DA$. Therefore the path above can be refined to a path of irreducible
morphisms. Convexity of $\zS $ in $\zG$ then implies that $X_i\in\zS_0$ for all
$i$. 
\qed
\end{pf}

\begin{cor}\label{cor 3.3}
Let  $\zS$ be a connected full subquiver in a connected component
$\zG$ of $\zG(\DA)$ such that $|\zS_0|= \textup{rk K}_0(A)$. The
following are equivalent:
\begin{itemize}
\item[(a)] $\zS$ is a section.
\item[(b)] $\zS$ is a local section.
\item[(c)] $\zS$ is a presection.
\end{itemize}
\end{cor}  
\begin{pf}
Again, it suffices to prove that (c) implies (a). Since $\zG$ is a
stable translation quiver, then the presection $\zS$ intersects each
$\tau_\DA$-orbit of $\zG$ at least once. Since, by hypothesis,
$|\zS_0|=\textup{rk K}_0(A)<\infty$, then $\zG$ is a transjective
component. The statement then follows from Proposition \ref{prop 2.4}.
\qed\end{pf}

\end{section}

\begin{section}{Local slices}\label{sect 4}
\begin{subsection}{Definition and examples}
We now define the main concept of this paper.

\begin{defn}\label{def local}
Let $C$ be a finite dimensional  algebra.
A \emph{local slice} $\zS$ in $\textup{mod}\, C$ is a local section
in a component $\zG$ of $ \zG(\textup{mod}\, C)$ such that 
 $|\zS| = \textup{rk K}_0(C)$.
\end{defn}

\begin{rem} 
Let $C$ be an algebra, and $\zS$ be a local slice in $\textup{mod}\,
C$.
\begin{itemize} 
\item[(a)] Since  local sections are presections, every path entirely
  contained in $\zS$ is sectional (because of Lemma \ref{lem 2.2}
  (b)). This implies that $\zS$ is acyclic.
\item[(b)] If $\zG$ is a stable component of $\zG(\textup{mod}\, C)$
  and $\zS\subset \zG$ then, by Lemma \ref{lem 2.2} (d), $\zS$
  intersects any $\tau_C$-orbit of $\zG$ at least once. Since
  $|\zS_0|<\infty$, this implies  that $\zG$ has only finitely many
  $\tau_C$-orbits. 
\end{itemize}   
\end{rem} 

\begin{exmp}\label{ex 3.3}
Let $C$ be a tilted algebra, and $\zG$ be a connecting component of
$\zG(\textup{mod}\, C)$. Then any complete slice in $\zG$ is also a
local slice. We shall prove below that any cluster-tilted algebra has
(many) local slices in its Auslander-Reiten quiver.
\end{exmp}

\begin{exmp}\label{ex 1}
The following is an example of an algebra which is neither tilted, nor
cluster-tilted, but whose Auslander-Reiten quiver contains a local slice.
Let $B$ be given by the quiver
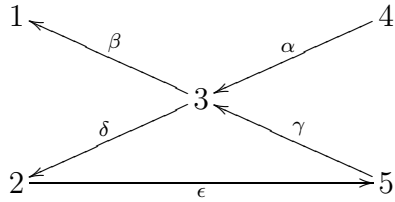
\[
\xymatrix@R=20pt@C=60pt{
1&&4\ar[ld]_\za\\
&3\ar[lu]_\zb\ar[ld]_\zd\\
2\ar[rr]_\ze&&5\ar[lu]_\zg
}
\]
bound by  $\za\,\zb=0$, $\zd\,\ze=0$, $\ze\,\zg=0$ and $\zg\,\zd=0$. 
The Auslander-Reiten quiver $\zG(\textup{mod}\,B)$ of $B$ is shown in
Figure \ref{fig arq}, 
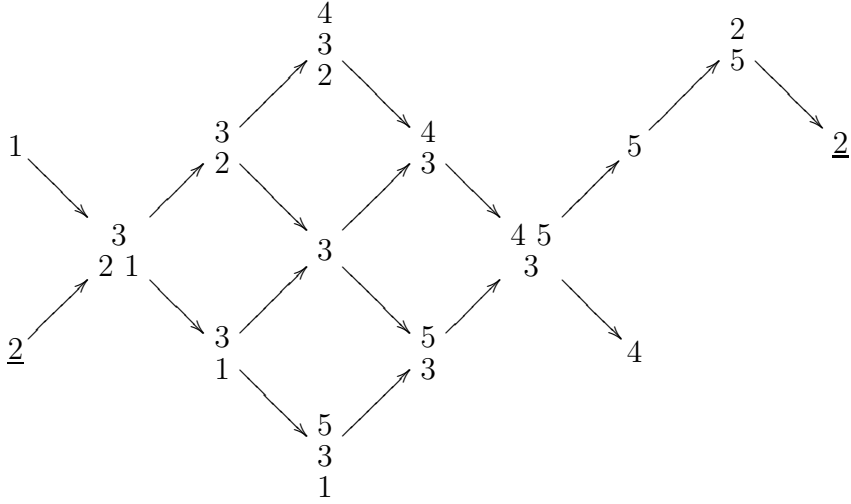
\begin{figure}
\[
\xymatrix@!R=15pt@!C=15pt{
&&& {\begin{array}{c}4 \vspace{-10pt} \\ \vspace{-10pt} 3\\2\end{array}} \ar[rd]&&&& {\begin{array}{c}2 \vspace{-10pt} \\ 5\end{array}} \ar[rd]\\
1\ar[rd]&&{\begin{array}{c}3 \vspace{-10pt} \\2\end{array}}\ar[ru]\ar[rd]
&&{\begin{array}{c}4 \vspace{-10pt} \\3\end{array}}\ar[rd]
&&\ 5\ \ar[ru]&& \ \underline{2}\ \\
&{\begin{array}{c}3 \vspace{-10pt} \\2\ 1\end{array}}\ar[ru]\ar[rd]
&&{\begin{array}{c}3\end{array}}\ar[ru]\ar[rd]
&&{\begin{array}{c}4\ 5 \vspace{-10pt} \\3\end{array}}\ar[ru]\ar[rd]\\
\underline{2}\ar[ru]&&{\begin{array}{c}3 \vspace{-10pt} \\1\end{array}}\ar[ru]\ar[rd]
&&{\begin{array}{c}5 \vspace{-10pt} \\3\end{array}}\ar[ru]&&4\\
&&&{\begin{array}{c}5 \vspace{-10pt} \\3 \vspace{-10pt} \\1\end{array}}\ar[ru]
}
\]
\caption{Auslander-Reiten quiver of Example \ref{ex 3.3}}\label{fig
  arq}
\end{figure}
where modules are represented by their Loewy series and we identify the
two copies of the underlined module $\underline{2}$.
Here $\zS=\left\{ \,{\begin{array}{c} \vspace{-10pt} 3\\ 2\ 1\end{array}} \,, \, {\begin{array}{c} \vspace{-10pt} 3\\ 2\end{array}} \,, \, 
{\begin{array}{c} \vspace{-10pt} 4\\ \vspace{-10pt} 3\\  2\end{array}} \,, \, {\begin{array}{c}3  \vspace{-10pt} \\ 1\end{array}} \,,  \,
{\begin{array}{c} \vspace{-10pt} 5\\ \vspace{-10pt} 3\\  1\end{array}} \,\right\}$ and 
 $\zS'=\left\{ \,{\begin{array}{c} 2 \vspace{-10pt} \\5  \end{array}} \,,\,
{\begin{array}{c}5\end{array}} \,,  \,
{\begin{array}{c}4\ 5 \vspace{-10pt} \\3  \end{array}} \,, \,
{\begin{array}{c}5 \vspace{-10pt} \\3 \end{array}} \,, \,
{\begin{array}{c}5 \vspace{-10pt} \\3 \vspace{-10pt} \\1 \end{array}} \,\right\}$
are local slices. 
Note that neither $\zS$ nor $\zS'$ is a section, because both
intersect twice the $\tau_B$-orbit of $\underline{2}$.
It is an interesting question to identify the algebras which have local slices.
\end{exmp}
\end{subsection}

\begin{subsection}{Cluster-tilted algebras}\label{sect 4.2}
Let $A$ be a  hereditary algebra. The \emph{cluster category}  $\CA$ of $A$
is defined as follows. Let $F$ be the automorphism of
$\DA=\DA^b(\textup{mod}\, A)$ defined as the composition
$\tau^{-1}_\DA [1]$, where $\tau^{-1}_\DA$ is the Auslander-Reiten
translation in $\DA$ and $[1]$ is the shift functor. Then $\CA$ is the
orbit category $\DA/F$, that is,  the objects of $\CA$ are the
$F$-orbits $\tilde X=(F^iX)_{i\in \mathbf{Z}}$, where $X\in\DA$, and
the set of morphisms from $\tilde X=(F^iX)_{i\in \mathbf{Z}}$ to
$\tilde Y=(F^iY)_{i\in \mathbf{Z}}$ is 
\[\Hom_{\CA}(\tilde X,\tilde Y) = \bigoplus_{i\in \mathbf{Z}}
\Hom_{\DA}(X,F^iY).
\]
It is shown in \cite{BMRRT,K} that $\CA$ is a triangulated
 category with almost split triangles. Furthermore, the projection functor $\pi_0:\DA\to \CA$ is a
 functor of triangulated categories and commutes with the
 Auslander-Reiten translations in both categories. We refer to \cite{BMRRT} for
 facts about the cluster category.

An object $\tilde T$ in $\CA$ is called a \emph{tilting object}
 provided $\Ext^1_{\CA}(\tilde T,\tilde T)=0$ and the number of
 isomorphism classes of indecomposable summands of $\tilde T$ equals
 $\textup{rk K}_0(A)$. The endomorphism algebra 
$B =\End_{\CA}(\tilde T)$ is then called a \emph{cluster-tilted
 algebra} \cite{BMR}. Of particular interest to us is the fact that
 the functor $\Hom_{\CA}(\tilde T, -):\CA\to \textup{mod}\, B$
 induces an equivalence 
\[\CA/\add(\tau \tilde T) \isomorphe \textup{mod}\, B,\]
where $\tau$ denotes the Auslander-Reiten translation in $\CA$, see
\cite{BMR}. This result entails several interesting consequences. For
instance, it is shown it \cite{KR} that any cluster-tilted algebra is
$1$-Gorenstein and hence of global dimension $1$ or $\infty$. For the
convenience of the reader, we give here a short proof of this fact.

\begin{prop}
Let $B$ be a cluster-tilted algebra.
\begin{itemize}
\item[(a)] For every injective $B$-module $I$, we have $\pd
  I\le 1$ (and for  every projective $B$-module $P$, we have
  $\textup{id}\, P\le 1$).
\item[(b)] $\textup{gl.dim.}\,B \in \{1,\infty\}$.
\end{itemize}   
\end{prop}  

\begin{pf} (a) By \cite[(IV.2.7) p.115]{ASS}, we need to prove that
  $\Hom_{B}(DB,\tau_{B}I)=0$. Now $\textup{mod}\,
  B\isomorphe \CA/\add(\tau \tilde T)$ (where $A$ and $T$ are
  as above) and, under this equivalence, every injective $\tilde
  C$-module is the image of an object of the form $\tau^2\tilde T_0\in
  \CA$, where $\tilde T_0\in \add\tilde T$.
It thus suffices to show that, for every  $\tilde T_0\in \add\tilde
  T$, we have $\Hom_{\CA}(\tau^2\tilde T,\tau^3\tilde T_0)=0$. But
  $\tau$ is an equivalence in $\CA$, hence the result follows from
  $\Hom_\CA(\tilde T,\tau\tilde T_0) \isomorphe \Ext^1_{\CA}(\tilde
  T,\tilde T_0)=0$.

(b) This is the proof of \cite{KR}, but we include it for
  completeness. It suffices to prove that, for every $B$-module $M$,
  id\,$M=d<\infty$ implies $\pd M\le 1$. Thus, let
\[
\xymatrix@C=20pt@R=0pt{
0\ar[r]&M\ar[r]&I^0\ar[r]^{f^1}&I^1\ar[r]&\cdots\ar[r]^{f^d}&I^d\ar[r]&0
}
\]
 be a minimal injective coresolution. Let $K^i=\textup{Im}\, f^i$ for every
 $i$. Then the exact sequence $0\to K^{d-1}\to I^{d-1}\to I^d\to 0$
 and (a) give $\pd K^{d-1}\le 1.$ An easy induction yields $\pd M\le 1$.
\qed
\end{pf}  

It is also shown in \cite{BMR} that the equivalence 
$\CA/\add(\tau \tilde T) \isomorphe \textup{mod}\, B$
commutes with the 
Auslander-Reiten translations in both categories.
Let  $\pi$ denote the composition of the functors 
\[
\xymatrix@C=90pt{\DA \ar@{->>}[r]^{\ \pi_0}&\ \CA\
  \ar@{->>}[r]^{\Hom_\CA(\tilde T,-)}& \textup{mod}\,B ,}\]
where $\pi_0$ is, as above, the canonical projection.
We notice that $\pi$ commutes with the 
Auslander-Reiten translations in both categories and also that, if
$X\in \DA$, then $\pi(X)=0$ if and only if $X\in\add(\tau_\DA\tilde
T)$.
\end{subsection}

\begin{subsection}{Auslander-Reiten quivers of cluster-tilted algebras}
With the above notations, we deduce the shape of the
Auslander-Reiten quiver of $\CA$ and  $\textup{mod}\, B$. 
 Let $Q$ be the ordinary quiver
of  $A$. If $A$ is representation-finite, then the Auslander-Reiten
quiver $\zG(\CA)$ is of the form $\mathbf{Z}Q/<\varphi>$, where
$\varphi$ is the automorphism of $\mathbf{Z}Q $ induced  by the
functor $F$. 
Since  $\zG(\CA)$ is stable and has sections isomorphic to $Q$, we say
that it is transjective.
If, on the other hand, $A$ is representation-infinite, then $\zG(\CA)$
consists of a unique component of the form $\mathbf{Z}Q$, which we
call transjective because it is the image under $\pi_0$ of the
transjective component of $\zG(\DA)$, and also of components which we
call regular because they are the image under $\pi_0$ of the
regular components of $\zG(\DA)$. In both cases, we deduce $\zG(\textup{mod}\, B)$ from
$\zG(\CA)$ by deleting the $|Q_0|$ points corresponding to the
summands of $\tau\tilde T$. In particular, $\zG(\textup{mod}\, \tilde
C)$ always has a unique transjective component, deduced from that of
$\zG(\CA)$ upon applying the functor $\Hom_\CA(\tilde T,-)$.

\begin{lem}\label{lem 3.3a} Let $\zG$ be a component of the
  Auslander-Reiten quiver of a cluster-tilted algebra $B$. If
  $\,\zG$ contains a local slice, then $\zG$ is the transjective component.
\end{lem}  

\begin{pf}
Assume that $\zG$ is not transjective. Then $\zG$ is either a stable
tube or a component of type $\mathbb{ZA}_\infty$, or is obtained from
one of these by deleting finitely many points. Also, since the functor
$\Hom_\CA(\tilde T,-):\CA\to \textup{mod}\, B$ commutes with
the Auslander-Reiten translations, deleting these points will not
change the $\tau$-orbits. Consequently, a local slice $\zS$  in $\zG$
lifts to a unique finite local section $\tilde \zS$ in a regular
component $\tilde \zG$ of $\zG(\CA)$. In particular, $\tilde \zG$ is
stable. Let thus $\tilde X\in\tilde \zS$ and $\tilde X=\tilde X_0\to
\tilde X_1\to\ldots\to \tilde X_i\to \ldots$ be  a sectional path of
irreducible morphisms (a ray) starting at $X$. Then $\tilde X_1$ or
$\tau \tilde X_1$ belongs to $\tilde \zS$. By induction, for each
$i\ge 0$, one of the objects $\tau^i\tilde X_i,\tau^{i-1}\tilde
X_{i},\ldots,\tilde X_i$ belongs to $\tilde \zS$. Therefore $\tilde
\zS$ is infinite, a contradiction.
\qed
\end{pf}  
\end{subsection}

\begin{subsection}{Lifting to the derived category}\label{sect 4.3}

\begin{lem}\label{lem 4.3} 
Let $ \zS$ be a  connected  full subquiver of the transjective
component of the Auslander-Reiten quiver of a cluster-tilted algebra  
$B$ arising from a hereditary algebra $A$ and  $\overline{\zS}$ be a
connected full subquiver of $D=\DA^b(\textup{mod}\,A)$ such that  
$\pi|_{\overline{\zS}}:\overline{\zS}\to \zS$ is bijective. Then 
$\zS$ is a local slice in $\textup{mod}\,B$ if and only if
$\overline{\zS}$ is a section in $\DA$ such that $|\zS_0|=\textup{rk
  K}_0(A)$. 
\end{lem}  

\begin{pf}
Since both subquivers are full, then the bijection
$\pi|_{\overline{\zS}}$ induces an isomorphism of quivers. Assume that
$\zS$ is a local slice in $\textup{mod}\,B$. We claim that
$\overline{\zS}$ is a presection in $\DA$.
 Assume 
that $\overline{X}\to \overline{Y}$ is an irreducible morphism in $\DA$ with 
$\overline{X}\in\overline{\zS_0} $.  Then we have two cases to
consider: 
\begin{enumerate}
\item If $\pi (\overline{Y})\ne 0$, then either $\pi
  (\overline{Y})\in\zS_0$ or $\tau_{B}\pi (\overline{Y})= \pi
  (\tau_\DA \overline{Y})\in
  \zS_0$. Therefore $\overline{Y}\in\overline{\zS_0}$ or
  $\tau_\DA\overline{Y}\in\overline{\zS_0}$.  
\item If $\pi (\overline{Y}) =0$, then $\pi(\tau_\DA \overline{Y}) \ne
  0$ because 
  $\Hom_{\DA}(\overline{Y},\tau_\DA \overline{Y}[1])\ne 0$.  
But we have an arrow $\overline{X}\to \overline{Y}$ which gives an
  arrow $\tau_\DA\overline{Y}\to 
  \overline{X}$,  hence an arrow $\pi(\tau_\DA\overline{Y})\to
  \pi(\overline{X})$. Since $\pi(\overline{X})\in \zS_0$, then
  $\pi(\tau_\DA \overline{Y})\in
  \zS_0$ and so $\tau_\DA\overline{Y}\in\overline{\zS}_0$. 
\end{enumerate}
Because of Lemma \ref{lem 2.2} (c), this shows that $\overline{\zS}$
is a presection. 
Since $\pi|_{\overline{\zS}}$ is a bijection, we have
$|\overline{\zS}_0|=|\zS_0|=\textup{rk K}_0(A)$. By Corollary
\ref{cor 3.3}, $\overline{\zS}$ is a section.

Conversely, assume that $\overline{\zS}$ is a section in $\DA$
and that $|\overline\zS_0|=\textup{rk K}_0(A)$. By Lemma \ref{lem
  3.2}, $\overline \zS$ lies in a transjective component of
$\zG(\DA)$. Note that
$|\overline\zS_0|=|\zS_0| $ and
$\textup{rk K}_0(A)=\textup{rk K}_0(B
)$.  Hence $|\zS_0|=\textup{rk K}_0(B)$.  

We show that $\zS $ is sectionally convex. 
Let $X=X_0\to X_1\to\ldots\to X_t=Y$ be a sectional path with $X,Y\in
\zS_0$. It lifts to a unique path 
$\overline{X}=\overline{X}_0\to \overline{X}_1\to\ldots\to
\overline{X}_t=\overline{Y}$  in $\DA$ with
$\overline{X},\overline{Y}\in \overline{\zS}_0$. Since $\overline{\zS}
$ is a section, then this path is sectional and all $\overline{X}_i\in
\overline{\zS}_0$. Applying $\pi$, we get that all $X_i$ lie in $\zS$.

Finally, we show that $\zS$ is a presection. Assume that
  $ X \to Y$ is an irreducible morphism in $\textup{mod}\, B$
  with $X\in \zS_0$. Let 
$\overline{X}=\pi|_{\overline{\zS}}^{-1}(X)$ and choose $\overline{Y}$ in
$\pi^{-1}(Y)$ such that we have an irreducible morphism
  $\overline{X}\to \overline{Y}$ in $\DA$. Then
  $\overline{Y}\in\overline{\zS}_0$ or 
$\tau_\DA\overline{Y}\in\overline{\zS}_0$. Hence $Y\in \zS_0$ or
  $\tau_{B} Y 
\in \zS_0$. Note that if $\tau_\DA \overline{Y} \in \overline \zS_0 $ then
$\pi(\tau_\DA\overline{Y})\ne 0$ because $\pi|_{\overline{\zS}} $ is a
bijection. We treat similarly the case of an irreducible morphism
  $X\to Y$ with $Y\in \zS_0$.
\qed\end{pf}

\end{subsection}

\begin{subsection}{Construction of local slices}\label{sect 4.4}

There is a close relation between tilted and cluster-tilted algebras. 
First, if $B$ is a cluster-tilted algebra, then there exist a
hereditary algebra $A$ and a tilting $A$-module $T$ such that 
$B\isomorphe \End_{\CA}(\tilde T)$, see
\cite[3.3]{BMRRT}. Also, if $A$ is a hereditary algebra and $T$ is a
tilting $A$-module so that the algebra 
 $C=\End_A(T)$ is tilted, the trivial extension 
 $\tilde C= C\ltimes \Ext^2_C(DC,C)$ (called the
\emph{relation-extension} of $C$) is cluster-tilted, 
and conversely, every cluster-tilted algebra is of this form, see
\cite{ABS}. 
Now,
since $\tilde C=C\ltimes \Ext^2_C(DC,C)$, then any $C$-module can be
considered as a $\tilde C$-module under the standard embedding 
 $i:\textup{mod}\,C\to \textup{mod}\,\tilde C$. Note that, in general,
$i$ does not preserve irreducible morphisms. 
We consider the complete slice $\zS=\add\Hom_A(T,DA)$ in
$\textup{mod}\, C$, where $C=\End \,T_A$ (see, for
instance, \cite{ASS} or \cite{R}) and denote its image as $i(\zS)=\zS'$ in
$\textup{mod}\, \tilde C$. The following Lemma collects the important
properties of $\zS'$.

\begin{lem}\label{lem 4.4}
Let  $\zS=\Hom_A(T,DA)$ and $\zS'=i(\zS)$, then
\begin{itemize}
\item[(a)] The image $\zS'$ is a local slice in $\textup{mod}\, \tilde C$.
\item[(b)] $i$ induces an isomorphism of quivers between $\zS$ and
  $\zS'$.
\item[(c)] $\Ann_{\tilde C}\,\zS' \cong\Ext^2_C(DC,C)$ as $C$-$C$-bimodules.
\end{itemize}
\end{lem}
\begin{pf}
(a) We have 
$\Hom_{\CA}(\tilde T,\widetilde{DA})=\Hom_{A}( T,DA)$ because $T$ is an
  $A$-module. Thus the image $i(\zS)=\zS'$ of $\zS$ is equal to the
  image  of $\zS$ under the composition  
\[
\xymatrix{
{\textup{mod}\,C} \ar[r]^{j\quad}& 
{\mathcal{D}^b(\textup{mod}\,C)}\ar[r]^{\phi}_{\cong}&
{\mathcal{D}^b(\textup{mod}\,A)}\ar[r]^{\quad \pi}&
{ \textup{mod}\,\tilde C,}
 }
\]
where $j$ is the inclusion in degree zero and $\phi$ the equivalence
$-\otimes_C^L T$. 
Let $\overline{\zS}$ be a connected full subquiver of
$\DA^b(\textup{mod}\, A)$ such that the restriction
$\pi|_{\overline{\zS}}:\overline{\zS}\to \zS'$ is bijective, that is,
  $\pi(\overline{\zS})=\zS'=i(\zS)$. Using the fact that $\phi$ is a
quasi-inverse of the derived functor 
  $R\textup{Hom}_A(T,-):\DA^b(\textup{mod}\, A)\to
  \DA^b(\textup{mod}\, C)$, we see that $\overline{\zS}=\phi\circ
j(\zS)$ is equal to $DA$ which is a connected section in a
transjective component of $\DA^b(\textup{mod}\, A)$. 
Since $\add(\tau\,T)\cap\add(DA)=\emptyset$, then  $\pi(I)\ne 0$ for
each $I\in \add DA$, so that $\pi|_{\add DA}:DA\to \zS'$ is bijective. 
 By Lemma \ref{lem 4.3}, $\zS'$ is a local slice in
 $\textup{mod}\,\tilde C$. 

$(b)$
This follows from the fact that $\zS'=\pi\phi j(\zS)$ and each one of
$j|_\zS,\phi$ and $\pi|_{\phi j(\zS)}$ preserves irreducible
morphisms (in the case of $\pi$, this is because $\add DA \cap
\add \tau\,T=\emptyset$). 

$(c)$ By \cite{ABS}, we have an isomorphism $\Ext^2_C(DC,C)
\isomorphe \Hom_{\DA}(T,FT)$ of
$C$-$C$-bimodules. Now 
\[
\begin{array}{rcl}
\Hom_A(T,DA)\cdot\Hom_{\DA}(T,FT)&=&\Hom_{\DA}(T,DA)\cdot\Hom_{\DA}(T,FT)\\  
&\subset& \Hom_{\DA}(T,F(DA)) \\ 
&=&\Hom_{\DA}(T,A[2])\ =\ 0,
\end{array}  \]
because $T$ is an $A$-module.
Therefore $\Hom_{\DA}(T,FT)\subset \Ann_{\tilde C} \,\zS'$.

Now let $f:T\to T$ be a non-zero morphism. Since its image is an
$A$-module, then $\Hom_A(\textup{Im} f,DA) \ne 0$. Let
$g:\textup{Im} f \to DA$ be a non-zero morphism and denote by $p:T\to
\textup{Im} f$ the canonical epimorphism.  Since $DA$ is an
injective module, there exists $g':T\to DA$ such that $g'f=gp\ne
0$. Therefore $\Ann_{\tilde C}\,\zS'\cap \Hom_A(T,T)=0$.  
Since, as a $k$-vector space, $\tilde C=\Hom_A(T,T)\oplus
\Hom_{\DA}(T,FT)$.
We have $\dim \tilde C = \dim\Hom_{\DA}(T,T)+\dim\Hom_{\DA}(T,FT) \le
\dim\Hom_{\DA}(T,T) +\dim \Ann_{\tilde C} \zS'$, because $\dim\Hom_{\DA}(T,FT)\subset
\Ann_{\tilde C}\zS'$. Since $\Hom_{\DA}(T,FT)$ and $\Ann_{\tilde C}
\zS'$ are in direct sum, then we have
$\dim\Hom_{\DA}(T,T) +\dim \Ann_{\tilde C} \zS'=\dim(\Hom_{\DA}(T,T)
+\Ann_{\tilde C} \zS')\le \dim
\tilde C$. Hence $\dim\Ann_{\tilde C}\zS'=\dim\Hom_{\DA}(T,FT)$. This 
 shows that
$\Hom_{\DA}(T,FT) $ and $\Ann_{\tilde C}\zS'$ are
equal as subspaces of $\tilde C$, and the statement follows.
\qed\end{pf}

\end{subsection} 

\begin{subsection}{Main result}\label{sect 4.5}

We are now ready for the proof of our main theorem.

\begin{thm}\label{thm 4.5}
Let $B$ be a cluster-tilted algebra. 
Then a tilted algebra $C$ is such that  $B=C\ltimes
\Ext^2_C(DC,C)$ if and only if there exists a local slice $\zS$ in
$\textup{mod}\, B$ such
that $C=B/\Ann_{B}\,\zS$. 
\end{thm}

\begin{pf} Necessity. 
It is well-known (see, for instance, \cite{ASS}) that any complete
slice $\zS$ in $\textup{mod}\, C$ is of the form
$\zS=\add\Hom_A(T,DA)$ for some hereditary algebra $A$ and some
tilting $A$-module $T$. 
By Lemma \ref{lem 4.4} (a), $\zS$
 embeds as a  local slice in
  $\textup{mod}\,B$ (because $B= \tilde C$). Moreover, by Lemma \ref{lem 4.4} (c), we have
  $\Ann_{B}\,\zS\isomorphe\Ext_C^2(DC,C)$ as $C$-$C$-bimodules.
 Therefore   
$C=B/\Ext_C^2(DC,C) =B/\Ann_{B}\,\zS$. 

Sufficiency. Let $B$ be cluster-tilted, and $\zS$ be a local slice in
$\textup{mod}\,B$. Set $C_1=B/\Ann_B\,\zS$.  Because of the definition
of a cluster-tilted algebra, there 
exist a hereditary algebra $A$ and a tilting object $\tilde T\in
\CA$ such that $B= \End_{\CA}(\tilde T)$. Let $\overline{\zS}$
be a connected component of the preimage  $\pi^{-1}(\zS)$ of $\zS$ in
$\DA^b(\textup{mod}\, A)$. Since the local slice $\zS$ can only occur in the
transjective component of $\zG(\textup{mod}\, B)$, because of
Lemma \ref{lem 3.3a}, then $\overline\zS$ lies in a transjective
component of $\zG(\DA^b(\textup{mod}\, A))$. Since $\zS$ and $\overline\zS$
have the same 
number of points, then $\pi|_{\overline \zS}:\overline\zS\to \zS$
is bijective, whence $\overline\zS$ is a section in
$\zG(\DA^b(\textup{mod}\, A))$ such that $|\overline\zS_0|=\textup{rk K}_0(A)$, by
Lemma \ref{lem 4.3}. We may
suppose without loss of generality that $\overline \zS =\add DA$. 

The fact that $\pi|_{\add DA}:\add DA\to \zS$ is a bijection implies
that the $F$-orbit  
 $\tau_\DA\,\pi^{-1}(\tilde T)$ in $\DA^b(\textup{mod}\, A)$ does not intersect
$\add DA$. Therefore, we have  $\add \pi^{-1} (\tilde T) \,\cap \,\add A[1]=\emptyset$, because
$\tau_\DA^{-1}DA=A[1]$ in $\DA^b(\textup{mod}\, A)$.
Thus, we can choose a representative $T$ in the $F$-orbit
$\pi^{-1}(\tilde T)$ such that $T$ is 
an $A$-module and $\pi(T)=\tilde T$. Then $T$ is a tilting $A$-module.
  Let $C_2=\End_A(T)$ be the corresponding tilted algebra.
  By \cite{ABS}, we have $B = 
C_2\ltimes \Ext^2_{C_2}(DC_2,C_2)$.
By Lemma \ref{lem 4.4} (c), we have  
$\Ext^2_{C_2}(DC_2,C_2)=\Ann_{B}\,\zS $.  Thus $C_2=B/\Ann_B\,\zS =C_1
$. In particular, $C_1$ is tilted and verifies $B =
C_1\ltimes \Ext^2_{C_1}(DC_1,C_1)$.
\qed
\end{pf}  

\begin{cor}\label{2.5}
Let $C$ be a tilted algebra and $\tilde C$ be the corresponding
cluster-tilted algebra. Then any complete slice in $\textup{mod}\, C$
embeds as a local slice in $\textup{mod}\, \tilde C$ and any local
slice in $\textup{mod}\, \tilde C$ arises this way.
\end{cor}  

\begin{pf}
This follows directly from Lemma \ref{lem 4.4} and Theorem \ref{thm 4.5}.
\qed
\end{pf}  
\end{subsection} 

\begin{subsection}{Computing the annihilator}\label{sect 4.7}
Our Theorem \ref{thm 4.5} actually gives a concrete way to compute the
tilted algebra $C$ starting from $\tilde C$. Given a local slice $\zS$
in $\textup{mod}\, \tilde C$, one computes its annihilator using the
following result.

\begin{cor}\label{cor 4.7}
Let $B$ be a cluster-tilted algebra and $\zS$ be a local slice
in $\textup{mod}\, B$. Then $\Ann_{B}\zS$ is generated
(as an ideal) by arrows in the quiver of $B$.
\end{cor}

\begin{pf} This follows from \cite[1.3]{AZ} using that $\tilde C$ is a
  trivial extension (hence a split extension) of $C$ by the
  $C$-$C$-bimodule $\Ext^2_C(DC,C)\isomorphe \Ann_{\tilde C}\zS. $
\qed
\end{pf}

We can be a bit more precise. Let $\tilde C=k\tilde Q/\tilde I$. Since
$\tilde C=\End_\CA\tilde T$, there is a bijection between the points
$x\in\tilde Q_0$ and the indecomposable summands $\tilde T_x$ of
$\tilde T$ so that each arrow $\za:x\to y $ in $ \tilde Q_1$
corresponds to a non-zero morphism 
$f_\za\in\Hom_\CA(\tilde T_y,\tilde
T_x)=\Hom_\DA(T_y,T_x)\oplus\Hom_\DA(T_y,FT_x)$. 
With this notation, $\Ann_{\tilde C}\zS$ is generated by all arrows
$\za:x\to y$ such that $f_\za\in \Hom_\DA(T_y,FT_x).$ This indeed
follows immediately from the isomorphisms 
\[\Hom_\DA(T,FT)\isomorphe
\Ext^2_C(DC,C) \isomorphe \Ann_{\tilde C} \zS.\]

Note that, as shown in \cite{ACT}, the arrows $\za$ which generate
$\Ann_{\tilde C}\zS$ have to satisfy certain
conditions. 
Moreover, if 
$C=kQ/I$, then $I=\tilde I\cap kQ$.
\end{subsection}

\begin{subsection}{Example}\label{sect 4.8}
Let $\tilde C$ be the cluster-tilted algebra (of type $\mathbb{D}_4$) given by
the quiver 
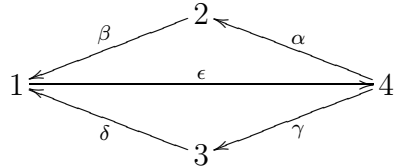
\[
\xymatrix@C=60pt@R=15pt{
&2\ar[dl]_\zb\\
1\ar[rr]^\ze&&4\ar[ul]_\za\ar[dl]^\zg\\
&3 \ar[ul]^\zd
}
\]
bound by $\za\zb=\zg\zd$,
 $\zb\ze=0$,
 $\zd\ze=0$,
 $\ze\za=0$,
 $\ze\zg=0$. 
The Auslander-Reiten quiver of $\tilde C$ is of the form shown in Figure \ref{fig arq2}
\begin{figure}
\[
\xymatrix@C=15pt@R=0pt{
&{\begin{array}{c} 2 \vspace{-10 pt} \\1 \end{array}}\ar[dr]
&&{\begin{array}{c} 3 \end{array}}\ar[dr]
&&{\begin{array}{c} 4 \vspace{-10 pt} \\2 \end{array}}\ar[dr]
&&{\framebox[20 pt]{$\tau\tilde T_2$}}
&&{\begin{array}{c} 2 \vspace{-10 pt} \\1 \end{array}}\\
{\begin{array}{c} \underline{1} \end{array}}\ar[dr]\ar[ur]
&{\framebox[20 pt]{$\tau\tilde T_4$}}
&{\begin{array}{c} 2\  3 \vspace{-10 pt} \\1 \end{array}}\ar[dr]\ar[ur]\ar[r]
&{\begin{array}{c} 4 \vspace{-10 pt} \\2\ 3 \vspace{-10 pt} \\1 \end{array}}\ar[r]
&{\begin{array}{c} 4 \vspace{-10 pt} \\2\ 3 \end{array}}\ar[dr]\ar[ur]
&{\framebox[20 pt]{$\tau\tilde T_1$}}
&{\begin{array}{c} 4 \end{array}}\ar[r]
&{\begin{array}{c} 1 \vspace{-10 pt} \\4 \end{array}}\ar[r]
&{\begin{array}{c} \underline{1} \end{array}}\ar[ru]\ar[rd]
& {\framebox[20 pt]{$\tau\tilde T_4$}}\\
&{\begin{array}{c} 3 \vspace{-10 pt} \\1 \end{array}}\ar[ru]&&
{\begin{array}{c} 2\end{array}}\ar[ru]&&
{\begin{array}{c} 4 \vspace{-10 pt} \\3 \end{array}}\ar[ru]&&
{\framebox[20 pt]{$\tau\tilde T_3$}}&&
{\begin{array}{c} 3 \vspace{-10 pt} \\1 \end{array}}
}
\]
\caption{Auslander-Reiten quiver of Example \ref{sect 4.8}}\label{fig arq2}
\end{figure}
where indecomposable modules are represented by their Loewy series,
and one identifies the two copies of  the underlined simple module
$1$ (and also the two copies of the modules ${\begin{array}{c} 2 \vspace{-10 pt} \\1
\end{array}}$ and ${\begin{array}{c} 3 \vspace{-10 pt} \\1 \end{array}}$). The entries
\framebox[20 pt]{$\tau\tilde T_i$} indicate the 
position of $\tau \tilde 
T_i$ in the cluster category. We find easily all the local slices
hence all the tilted 
subalgebras of $\tilde C$ which realise it as a
relation-extension. There are only the three algebras $C_i$ shown in 
Figure \ref{fig arq3} each corresponding to a local slice $\zS_i$.
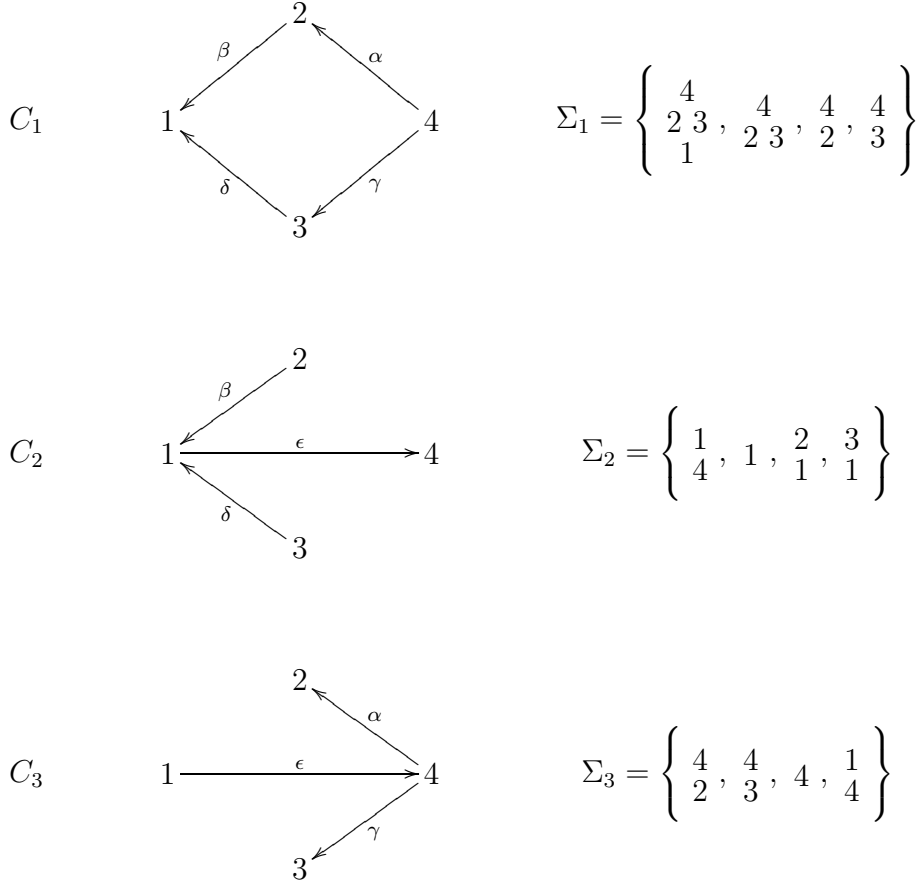
\begin{figure}
\[
\xymatrix@C=40pt@R=10pt{
&&2\ar[dl]_\zb\\
C_1&1&&4\ar[ul]_\za\ar[dl]^\zg
&{\zS_1=\left\{ \,\begin{array}{c}4 \vspace{-10 pt} \\2\ 3 \vspace{-10 pt} \\1\end{array}  \, , \,
\begin{array}{c}4 \vspace{-10 pt} \\2\ 3\end{array}  \,,  \,
\begin{array}{c}4 \vspace{-10 pt} \\2\end{array}  \,,  \,
\begin{array}{c}4 \vspace{-10 pt} \\3\end{array} \, \right\}}\\
&&3 \ar[ul]^\zd\\
\\
\\
&&2\ar[dl]_\zb\\
C_2&1\ar[rr]^\ze&&4
&{\zS_2=\left\{ \,\begin{array}{c} \vspace{-10 pt} 1\\4\end{array}  \,,  \,
\begin{array}{c}1\end{array} \, ,  \,
\begin{array}{c}2 \vspace{-10 pt} \\1\end{array}  \,,  \,
\begin{array}{c}3 \vspace{-10 pt} \\1\end{array} \, \right\}}\\
&&3 \ar[ul]^\zd
\\
\\
\\
&&2\\
C_3&1\ar[rr]^\ze&&4\ar[ul]_\za\ar[dl]^\zg
&{\zS_3=\left\{ \,\begin{array}{c}4 \vspace{-10 pt} \\2\end{array} \, ,  \,
\begin{array}{c}4 \vspace{-10 pt} \\3\end{array}  \,,  \,
\begin{array}{c}4\end{array} \, ,  \,
\begin{array}{c}1 \vspace{-10 pt} \\4\end{array} \, \right\}}\\
&&3
}
\]
\caption{Tilted algebras and corresponding local slices}\label{fig arq3}
\end{figure}
 with the inherited relations in each case.

\end{subsection} 
\end{section}

\begin{section}{Cluster-tilted algebras of tree type}\label{sect 5}
In view of our main result, Theorem \ref{thm 4.5}, it is reasonable
to ask whether there exist sufficiently many local slices in the
Auslander-Reiten quiver of a cluster-tilted algebra. Since the latter
is deduced from the Auslander-Reiten quiver of the cluster category by
dropping finitely many points, and since local slices can only occur
in the transjective component, then all but at most finitely many
indecomposables lying in the transjective component of the
cluster-tilted algebra belong to a local slice. That not necessarily
all indecomposables in the transjective component belong to a local
slice is seen in the following example.

\begin{exmp} 
Let $A$ be the hereditary algebra given by the quiver

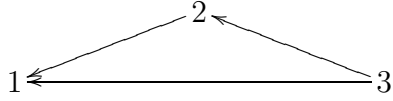
\[
\xymatrix@C=60pt@R=15pt{
&2\ar[ld]\\
1&&3\ar[lu]\ar[ll]
}
\]
 and consider the tilting $A$-module $T= \ 1\ \oplus 
\begin{array}{c}3 \vspace{-10 pt} \\1\ 2 \vspace{-10 pt} \\ \ \ \ 1 
\end{array}  
\oplus 
\begin{array}{c}
3 \vspace{-10 pt} \\1
\end{array}  $.
Note that while $T_1=1$ and $T_2=
\begin{array}{c} 3 \vspace{-10 pt} \\1\ 2 \vspace{-10 pt} \\\ \ \ 1 
\end{array}  $
are projective $A$-modules, the module $T_3=
\begin{array}{c}3 \vspace{-10 pt} \\1
\end{array}$
is simple regular non-homogeneous. The transjective component of the
Auslander-Reiten quiver of $\tilde C=\End_\CA \tilde T$ is thus of the form

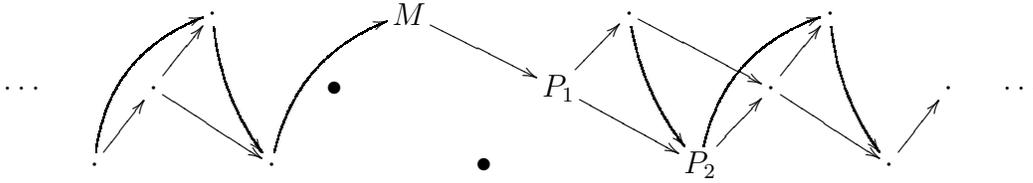
\[ 
\xymatrix@R=15pt@C=15pt{
&&&\cdot\ar@/_5pt/[rdd]&&&M\ar[rrd]&&&\cdot\ar[rrd]\ar@/_5pt/[rdd]&&&\cdot\ar@/_5pt/[rdd]&\\
\ldots&&\cdot\ar[ru]\ar[rrd]&&&\bullet&&&P_1\ar[ru]\ar[rrd]&&&\cdot\ar[ru]\ar[rrd]&&&\cdot&\ldots\\
&\cdot\ar[ru]\ar@/^12pt/[rruu]&&&\cdot\ar@/^12pt/[rruu]&&&\bullet&&&P_2\ar[ru]\ar@/^12pt/[rruu]&&&\cdot\ar[ru]
}
\]
where $P_i=\Hom_\CA(\tilde T,T_i)$, for  $i\in \{1,2\}$, and the two
$\bullet$ represent $\tau T_1$ and $\tau T_2$. It is easily seen that
the indecomposable $\tilde C$-module $M=\rad P_1$ lies on no local slice.

\end{exmp}

We say that a cluster-tilted algebra $\End_{\CA}(\tilde T)$ is of
\emph{tree type} if the ordinary quiver of $A$ is a  tree.

\begin{thm}\label{2.6}
Let $\tilde C$ be a cluster-tilted algebra of tree type. Then any
indecomposable $\tilde C$-module lying in the transjective component
lies on a local slice.
\end{thm}   

\begin{pf}
Let $\tilde M$ be an indecomposable $\tilde C$-module in the
transjective component and   $M$ be an indecomposable object in  
$\pi^{-1}(\tilde M)$; here, as in section \ref{sect 4}, $\pi$ denotes
the composition of the natural functors 
\[\xymatrix@C=80pt{\DA^b(\textup{mod}\, A)\ar@{>>}[r]^{\pi_0}&\CA\ar@{>>}[r]^{\Hom_\CA(\tilde T,-)}&\textup{mod}\,
  \tilde C}.\]
We claim that $M$ lies on  a section $\zS$ in $\DA^b(\textup{mod}\,
A)$
such that $|\zS_0|=\textup{rk K}_0(A)$ and 
$\zS\cap\add\tau\,\tilde T=\emptyset$.

Before proving this claim, we show that the
theorem follows from it. Indeed, under these conditions,
 $\pi|_{\zS}:\zS\to\pi(\zS)$ is 
bijective and $\pi(\zS)$ is a  connected full subquiver of the
transjective component of the Auslander-Reiten quiver of $\tilde
C$. By Lemma \ref{lem 4.3}, $\pi(\zS)$ is a local
slice. Obviously $\tilde M=\pi(M)\in\pi(\zS)$.

We now prove the claim.
First, we fix some terminology.
For any connected full subquiver $\mathcal{T}$ of
$\zG(\DA^b(\textup{mod}\, A))$ such that
$\mathcal{T}$ is a tree and $M\in \mathcal{T}_0$, there is a
unique reduced walk in $\mathcal{T}$ from $M$ to any other point
$N\in\mathcal{T}_0$. We define the distance $d_{\mathcal{T}}(M,N)$ between $M$
and $N$ in $\mathcal{T}$ to be the number of arrows of this walk.
Moreover, we define  
\[ d_{\mathcal{T}} =  \left\{
\begin{array}{ll} \textup{rk K}_0 (A) & \textup{if $\mathcal{T}\cap
  \add \tau\,\tilde T=\emptyset$}\\
\min\{d_{\mathcal{T}}(M,N)\mid N\in\mathcal{T}\cap \add \tau\,\tilde T\}
  & \textup{if $\mathcal{T}\cap  \add \tau\,\tilde T\ne\emptyset$.}
\end{array}\right.  
\]
Thus $d_{\mathcal T}$ measures the distance between $M$ and the set
$\add\tau\,\tilde T$ in the tree $\mathcal{T}$.
Note that $d_{\mathcal{T}}\ge 1$ because $M\notin \add\tau\,\tilde T$ and
$d_{\mathcal{T}}= \textup{rk K}_0 (A)$ if and only if $\mathcal{T}\cap
\add  \tau\,\tilde T=\emptyset$.

Now, let $\zS_1$ be any section in $\DA^b(\textup{mod}\, A)$ containing
$M$. By hypothesis, $\zS_1$ is a tree. 
If $d_{\zS_1}=  \textup{rk K}_0 (A)$ we are done. Suppose that
$d_{\zS_1}<  \textup{rk K}_0 (A)$.
Consider the set  $\mathcal{N}=\{N\in\zS_1\cap\add\tau\,\tilde T \mid
d_{\zS_1}=d_{\zS_1}(M,N)\}$ and let $N\in\mathcal{N}$. 
Consider the unique reduced walk 
\[M=X_{0}\frac{\quad}{\quad} X_{1} \frac{\quad}{\quad} \ldots 
\frac{\quad}{\quad}
 X_{(d_{\zS_1}-1)}=L \frac{\quad}{\quad}X_{d_{\zS_1}}=N \]
 in $\zS_1$ from $M$ to $N$.
Deleting the edge  $L \frac{\quad}{\quad} N$ cuts $\zS_1$ into two subtrees.
Let $\zS_1^M$ be the subtree containing $M$ (and $L$) and let  $\zS_1^N$ the
 subtree containing $N$.
There are two cases to consider, according to the orientation of the
 arrow between $L$ and $N$.
\begin{enumerate}
\item If $L\to N$, we define $\zS_2$
to be the full subquiver of $\zG(\DA^b(\textup{mod}\, A))$ having as
points those of $\zS_1^M\cup \tau\,\zS_1^N$.
\item If $L\ot N$, we define  $\zS_2$
to be the full subquiver of $\zG(\DA^b(\textup{mod}\, A))$ having as
points those of
$\zS_1^M\cup \tau^{-1}\,\zS_1^N$.  
\end{enumerate}   
By construction, $\zS_2$ is a connected tree, it lies in the transjective
component and it intersects every 
$\tau$-orbit exactly once.

We now show that $\zS_2$ is also convex. 
Assume that $L\to N$ (the proof in case $L\ot N$ is entirely similar).
 Then $\zS_2$ has two subtrees $\zS_1^M$ and
$\tau\,\zS_1^N$, connected by the arrow $
\tau\,N\to L$. We show first that these two subtrees are convex.
Suppose  that $X=X_0\to
X_1\to\ldots\to X_s=Y$ is a path in $\zG(\DA^b(\textup{mod}\, A)$ with
$X,Y\in\zS_1^M$.
By convexity of $\zS_1$, we have $X_i\in\zS_1$ for all $i$.
Since there is exactly  one  walk from $X$ to $Y$ in $\zS_1$ we actually
have $X_i$ in $\zS_1^M$ for all $i$. Thus $\zS_1^M$ is
convex. Similarly,  $\tau\,\zS_1^N$ is also convex.
Now suppose  that $X=X_0\to
X_1\to\ldots\to X_s=Y$ is a path with $X,Y\in\zS_2$.
If $X\in\zS_1^M$, then $Y\in \zS_1^M$ because of the structure of the
transjective component, and then the convexity of $\zS_1^M$  implies
the result.
If $X,Y \in\tau\,\zS_1^N$ we are done by convexity of
$\tau\,\zS_1^N$. Hence the only remaining case is when
$X\in\tau\,\zS_1^N$ and $Y\in\zS_1^M$.
Suppose that there is an $i\in\{1,\ldots,s-1\}$ such that  $X_i\notin
\zS_1$. We may suppose, without loss of generality, that $i=1$.
Now $\tau\,\zS_1$ is a section such that $|(\tau\,\zS_1)_0|=\textup{rk
  K}_0(A)$, hence, in particular, it is a presection. Therefore
$X\in\tau\,\zS_1, \ X\to X_1$ and 
$X_1\notin\tau\,\zS_1$ imply $\tau\,X_1\in\tau\,\zS_1$, hence
$X_1\in\zS_1$. Then all $X_i$ with $1\le i\le s$ lie in $\zS_1$
because $\zS_1$ is convex. Since $X_1\notin\zS_2 $, we have
$X_1\in\zS_1^N$. Then, since $Y\in\zS_1^M$ and the subtrees $\zS_1^N,\zS_1^M$
are only joined by the arrow $L\to N$, there exists $j$ such that
$X_j=N$ and $X_{j+1}=L$. But then there is an arrow $L\ot N$,  a contradiction.
 This shows that $\zS_2$ is convex and hence is a section in
 $\DA^b(\textup{mod}\, A)$
satisfying $|(\zS_2)_0|=\textup{rk K}_0(A).$

We repeat this construction for every element $N$ of
$\mathcal{N}$. Note 
that, since $N\in\add \tau\,\tilde T$, neither $\tau 
\,N$ nor $\tau^{-1}N$ are in  $\add \tau\,\tilde T$.
 In this way, we obtain, after $k=|\mathcal{N}|$ steps,
a section $\zS_{k+1}$  in $\zG(\DA^b(\textup{mod}\, A))$ that contains
$M$, such that $|(\zS_{k+1})_0|=\textup{rk K}_0(A)$ and 
\[d_{\zS_{k+1}}>d_{\zS_1}.\]
If $d_{\zS_{k+1}}=\textup{rk K}_0 (A)$ then $\zS_{k+1}\cap
\add(\tau\,\tilde T) =\emptyset$ and we are done.
Otherwise, we repeat the construction with the set 
\[\mathcal{N}_{k+1}=\{N\in\zS_{k+1}\mid d_{\zS_{k+1}}(M,N)=d_{\zS_{k+1}}\}.\]
This algorithm will find the required
section $\zS$ in a finite number of steps.
\qed
\end{pf}  
 \begin{cor}\label{cor 5.3}
Let $\tilde C$ be a representation-finite cluster-tilted algebra. Then
any indecomposable $\tilde C$-module lies on a  local slice.
\end{cor}

\end{section}

\end{document}